\documentclass[a4paper,10pt]{article}
\usepackage[utf8x]{inputenc}
\usepackage{amsthm,amsmath,amssymb,oldgerm}
\usepackage[margin=3cm]{geometry}
\theoremstyle{plain}
\newtheorem{thm}{Theorem}[section]
\newtheorem{lem}[thm]{Lemma}
\newtheorem{prop}[thm]{Proposition}
\newtheorem{cor}[thm]{Corollary}
\theoremstyle{definition}
 
\theoremstyle{definition}
\newtheorem{rem}[thm]{Remark}

\let\Im\relax
\DeclareMathOperator{\Im}{Im}
\DeclareMathOperator{\del}{\Delta_{hyp}}
\DeclareMathOperator{\hyp}{\mu_{hyp}} 
\DeclareMathOperator{\can}{\mu_{can}}
\DeclareMathOperator{\hathyp}{\widehat{\mu}_{hyp}}
\DeclareMathOperator{\hatcan}{\widehat{\mu}_{can}}
\DeclareMathOperator{\shyp}{\mu_{shyp}}
\DeclareMathOperator{\hatshyp}{\widehat{\mu}_{shyp}}

\DeclareMathOperator{\khyp}{{\it{K_{\mathrm{hyp}}}}}

\DeclareMathOperator{\ghyp}{\it{g_{\mathrm{hyp}}}}
\DeclareMathOperator{\gcan}{\it{g_{\mathrm{can}}}}

\let\Re\relax
\DeclareMathOperator{\Re}{Re}
\DeclareMathOperator{\vx}{\mathrm{vol_{\mathrm{hyp}}}}

\DeclareMathAlphabet{\mathpzc}{OT1}{pzc}{m}{it}
\makeatletter

\newcommand{\Rmnum}[1]{\expandafter\@slowromancap\romannumeral #1@}
\title{Extension of a key identity}
{\small\author{Anilatmaja Aryasomayajula}}
\begin{document}
\maketitle
\begin{abstract}
\noindent In this article, we extend a certain key identity proved by J.~Jorgenson and J.~Kramer 
in \cite{jk} to noncompact hyperbolic Riemann orbisurfaces of finite volume. This identity relates the two 
natural metrics, namely the hyperbolic metric and the canonical metric defined on a Riemann orbisurface. 
\end{abstract}
\section*{Introduction}
\paragraph{Notation}
Let $X$ be a noncompact hyperbolic Riemann orbisurface of finite volume $\vx(X)$ with genus $g\geq 1$, and 
can be realized as the quotient space $\Gamma\backslash\mathbb{H}$, where $\Gamma\subset \mathrm{PSL}_{2}(
\mathbb{R})$ is a Fuchsian subgroup of the first kind acting on the hyperbolic upper half-plane $\mathbb{H}$, 
via fractional linear transformations. 

Let $\mathcal{P}$ denote the set of cusps of $\Gamma$, and put $\overline{X}=X\cup \mathcal{P}$. Then, 
$\overline{X}$ admits the structure of a Riemann surface. 

\vspace{0.2cm}
Let $\hyp$ denote the (1,1)-form associated to hyperbolic metric, which is 
the natural metric on $X$, and of constant negative curvature minus one. 

\vspace{0.2cm}
The Riemann surface $\overline{X}$ is embedded in its Jacobian variety $\mathrm{Jac}(\overline{X})$ via the 
Abel-Jacobi map. Then, the pull back of the flat Euclidean metric by the Abel-Jacobi map is called the 
canonical metric, and the (1,1)-form associated to it is denoted by $\hatcan$. We denote its restriction 
to $X$ by $\can$.

\vspace{0.2cm}
Let $\del$ denote the hyperbolic Laplacian acting on smooth functions on $X$. Let $\khyp(t;z,w)$ denote the 
hyperbolic heat kernel defined on $\mathbb{R}_{>0}\times X\times X$, which is the unique solution of the 
heat equation 
\begin{align*}
\bigg(\Delta_{\text{hyp},z} + \frac{\partial}{\partial t}\bigg)\khyp(t;z,w) =0,
\end{align*}
and the normalization condition
\begin{align*}  
\lim_{t\rightarrow 0}\int_{X}\khyp(t;z,w)f(z)\hyp(z) = f(w),
\end{align*}
for any fixed $w\in X$ and any smooth function $f$ on $X$. When $z=w$, for brevity of notation, we denote the 
hyperbolic heat kernel by $\khyp(t;z)$.

\vspace{0.2cm}
Let $C_{\ell,\ell\ell}(\overline{X})$ denote the space of singular functions, which are $\log$-singular at 
finitely many points of $X$, and are $\log\log$-singular at the cusps. With notation as above, we now state the 
main result.
\paragraph{Main result}
With notation as above, for any $f\in C_{\ell,\ell\ell}(\overline{X})$, we have the equality of integrals 
\begin{align*}
&g\int_{X}f(z)\can(z) =\notag \\&\bigg(\frac{1}{4\pi}+\frac{1}{\vx(X)}\bigg)\int_{X}f(z)\hyp(z) + 
\frac{1}{2}\int_{X}f(z)\bigg(\int_{0}^{\infty}\del \khyp(t;z)dt\bigg)\hyp(z).
\end{align*}
The above relation, which relates the two natural metrics defined on a Riemann orbisurface has been 
proved for compact hyperbolic Riemann surfaces, as a relation of differential forms by J.~Jorgenson and 
J.~Kramer in \cite{jk}. The same authors have also extended the key identity to noncompact hyperbolic Riemann 
surfaces of finite volume in \cite{bonn}. In this paper, the authors use different methods from \cite{jk}, and study the 
behavior of the key identity over a family of degenerating compact hyperbolic Riemann surfaces. 

\vspace{0.2cm}
Our main theorem can be seen as an extension of their result to 
elliptic fixed points and cusps at the level of currents acting on the space of singular function 
$C_{\ell,\ell\ell}(\overline{X})$. Our methods are different from the ones employed in \cite{bonn}, and are organized 
around the original line of proof in \cite{jk}.
\paragraph{Arithmetic significance}
The key identity has been the most significant technical result of \cite{jk}, which transforms a 
problem in Arakelov theory into that of hyperbolic geometry. The key identity has enabled J.~Jorgenson and 
J.~Kramer to derive optimal bounds for the canonical Green's function defined on a compact hyperbolic Riemann 
surface $X$ in terms of invariants coming from the hyperbolic geometry of $X$. These 
bounds were essential for B.~Edixhoven's algorithm in \cite{edix} for computing certain Galois 
representations associated to a fixed modular form of arbitrary weight. 

Furthermore, using the key identity and the Polyakov formula, J.~Jorgenson and J.~Kramer have obtained 
optimal bounds for the Faltings delta function in \cite{jkannals}. The key identity is again the most 
important technical tool. 

\vspace{0.2cm}\noindent
Using the key identity one can relate the holomorphic world of cusp forms with the $C^{\infty}$ 
world of M\"ass forms, via the spectral expansion of the hyperbolic heat kernel $\khyp(t;z)$ in terms 
of M\"ass forms. In fact, J.~Jorgenson and J.~Kramer have derived a Rankin-Selberg $L$-function 
relation relating the Fourier coefficients of cusp forms with those of M\"ass forms in \cite{bonn}.

\vspace{0.2cm}\noindent
The extended version of the key identity enables us to extend the work of J.~Jorgenson and J.~Kramer to 
noncompact hyperbolic Riemann orbisurfaces of finite volume. In an upcoming article \cite{anil}, using the 
key identity, we extend the bounds derived in \cite{jk}.
{\small{\paragraph{Acknowledgements}
This article is part of the PhD thesis of the author, which was completed under the supervision of  
J.~Kramer at Humboldt Universit\"at zu Berlin. The author would like to express his gratitude to J.~Kramer 
for his support and many valuable scientific discussions. The author would also like to extend his 
gratitude to J.~Jorgenson for sharing new scientific ideas, and to R.~S.~de~Jong for many interesting 
scientific discussions and for pointing out a mistake in the first proof.}}  
\section{Background material}\label{section1}
Let $\Gamma \subset \mathrm{PSL}_{2}(\mathbb{R})$ be a Fuchsian subgroup of the first kind acting by 
fractional linear transformations on the upper half-plane $\mathbb{H}$. Let $X$ be the quotient space 
$\Gamma\backslash \mathbb{H}$, and let $g$ denote the genus of $X$. The quotient space $X$ admits the 
structure of a Riemann orbisurface.  

Let $\mathcal{E}$, $\mathcal{P}$ be the finite set of elliptic fixed points and cusps of $X$, respectively; put $\mathcal{S}
 = \mathcal{E}\cup\mathcal{P}$. For $\mathfrak{e}\in \mathcal{E}$, let $m_{\mathfrak{e}}$ denote the order of $\mathfrak{e}$; for 
$p\in \mathcal{P}$, put $m_{p}=\infty$; for $z\in X\backslash \mathcal{E}$, put $m_{z}=1.$ Let $\overline{X}$ denote $\overline{X}=X\cup
\mathcal{P}.$

Locally, away from the elliptic fixed points and cusps, we identity $\overline{X}$ with its universal cover 
$\mathbb{H}$, and hence, denote the points on $\overline{X}\backslash \mathcal{S}$ by the same letter as the points on $\mathbb{H}$.
\paragraph{Structure of $\overline{X}$ as a Riemann surface} 
The quotient space $\overline{X}$ admits the structure of a compact Riemann surface. We refer the reader to 
section 1.8 in $\cite{miyake}$, for the details regarding the structure of $\overline{X}$ as a compact 
Riemann surface. For the convenience of the reader, we recall the coordinate functions for the neighborhoods 
of elliptic fixed points and cusps.

Let $w\in U_{r}(\mathfrak{e})$ denote a coordinate disk of radius $r$ around an elliptic fixed point $\mathfrak{e}\in\mathcal{E}$. Then, the coordinate function 
$\vartheta_{\mathfrak{e}}(w)$ for the coordinate disk $ U_{r}(\mathfrak{e})$ is given by
\begin{equation*}\label{ellipticcoordinatefn}
\vartheta_{\mathfrak{e}}(w)= \bigg(\frac{w-\mathfrak{e}}{w-{\overline{\mathfrak{e}}}}\bigg)^{m_{\mathfrak{e}}}.
\end{equation*}
Similarly, let $p\in\mathcal{P}$ be a cusp and let $w\in U_{r}(p)$. Then $\vartheta_{p}(w)$ is given by
\begin{equation*}
 \vartheta_{p}(w)= e^{2\pi i \sigma_{p}^{-1}w},
\end{equation*}
where $\sigma_{p}$ is a scaling matrix of the cusp $p$ satisfying the following relations
\begin{align*}
\sigma_{p}i\infty = p \quad \mathrm{and} \quad \sigma_{p}^{-1}\Gamma_{p}\sigma_{p} = \langle\gamma_{\infty}\rangle,\quad
\mathrm{where}\,\,\, \gamma_{\infty}=\left(\begin{array}{ccc} 1 & 1\\ 0 & 1  \end{array}\right)
\quad&\mathrm{and}\quad\Gamma_{p}=\langle\gamma_{p}\rangle
\end{align*}
denotes the stabilizer of the cusp $p$ with generator $\gamma_{p}$.
\paragraph{Hyperbolic metric} 
We denote the (1,1)-form corresponding to the hyperbolic metric of $X$, which is compatible with the complex 
structure on $X$ and has constant negative curvature equal to minus one, by $\hyp(z)$. Locally, for 
$z\in X\backslash \mathcal{E}$, it is given by
\begin{equation*}
 \hyp(z)= \frac{i}{2}\cdot\frac{dz\wedge d\overline{z}}{{\Im(z)}^{2}}.
\end{equation*} 
Let $\vx(X)$ be the volume of $X$ with respect to the hyperbolic metric $\hyp$. It is given by the formula 
\begin{equation*}
\vx(X) = 2\pi\bigg(2g -2 + |\mathcal{P}|+\sum_{\mathfrak{e} \in \mathcal{E}}\bigg(1-
\frac{1}{m_{\mathfrak{e}}}\bigg)\bigg). 
\end{equation*}
The hyperbolic metric $\hyp(z)$ is singular at the elliptic fixed points and at the cusps, and defines a 
singular and integrable (1,1)-form on $\overline{X}$, which we denote by $\hathyp(z)$. The rescaled 
hyperbolic metric 
\begin{equation*}
 \shyp(z)= \frac{\hyp(z)}{ \vx(X)}
\end{equation*}
measures the volume of $X$ to be one, and we denote the (1,1)-form determined by $\shyp(z)$ on $\overline{X}$ 
by $\hatshyp(z)$. Furthermore, let us denote the (1,1)-currents determined by $\hathyp(z)$ and $\hatshyp(z)$ acting on smooth functions defined on 
$\overline{X}$ by $[\hathyp(z)]$ and $[\hatshyp(z)]$, respectively. 

Locally, for $z$ $\in$ $X$, the hyperbolic Laplacian $\Delta_{\mathrm{hyp}}$ on  $X$ is given by
\begin{equation*}
\Delta_{\mathrm{hyp}} = -y^{2}\bigg(\frac{\partial^{2}}{\partial x^{2}} +
 \frac{\partial^{2}}{\partial y^{2}}\bigg) = -4y^{2}\bigg(\frac{\partial^{2}}{\partial z
\partial \overline{z}} \bigg). 
\end{equation*}
Recall that $d=\left(\partial + \overline{\partial} \right), 
$ $d^{c}=\dfrac{1}{4\pi i}\left( \partial - \overline{\partial}
\right)$, and $dd^{c}= -\dfrac{\partial\overline{\partial}}{2\pi i}$. 
\paragraph{Canonical metric}
Let $S_{2}(\Gamma)$ denote the $\mathbb{C}$-vector space of cusp forms of 
weight 2 with respect to $\Gamma$ equipped with the Petersson inner product. Let 
$\lbrace f_{1},\ldots,f_{g}\rbrace $ denote an orthonormal basis of $S_{2}(\Gamma)$ with respect to the Petersson 
inner product. Then, the (1,1)-form $\can(z)$ corresponding to the 
canonical metric of $X$ is given by 
\begin{equation*}
 \can(z)=\frac{i}{2g} \sum_{j=1}^{g}\left|f_{j}(z)\right|^{2}dz\wedge d\overline{z}.
\end{equation*}
The canonical metric $\can(z)$ remains smooth at the elliptic fixed points and at the cusps, and measures 
the volume of $X$ to be one. We denote the smooth (1,1)-form defined by $\can(z)$ on $\overline{X}$ by $\hatcan(z)$, and the (1,1)-current 
determined by $\hatcan(z)$ acting on smooth functions defined on $\overline{X}$ by $[\hatcan(z)]$. 
\paragraph{Canonical Green's function}
For $z, w \in \overline{X}$, the canonical Green's function $\widehat{g}_{\mathrm{can}}(z,w)$ is defined as 
the solution of the differential equation 
\begin{equation}\label{diffeqngcan}
d_{z}d^{c}_{z}\widehat{g}_{\mathrm{can}}(z,w)+ \delta_{w}(z)=\hatcan(z),
\end{equation}
with the normalization condition
\begin{equation*}\label{normcondgcan}
 \int_{\overline{X}}\widehat{g}_{\mathrm{can}}(z,w)\hatcan(z)=0. 
\end{equation*}
From equation \eqref{diffeqngcan}, it follows that $\widehat{g}_{\mathrm{can}}(z,w)$ admits a 
$\log$-singularity at $z=w$, i.e., for $z, w\in \overline{X}$, it satisfies 
\begin{equation}\label{gcanbounded}
\lim_{w\rightarrow z}\big( \widehat{g}_{\mathrm{can}}(z,w)+ 
\log |\vartheta_{z}(w)|^{2}\big)= O_{z}(1). 
\end{equation}
For a fixed $w\in\overline{X}$, the canonical Green's function $\widehat{g}_{\mathrm{can}}(z,w)$ determines a 
current $[\widehat{g}_{\mathrm{can}}(\cdot,w)]$ of type (0,0) acting on smooth (1,1)-forms defined on 
$\overline{X}$. Furthermore, for a fixed $w\in\overline{X}$, the current $[\widehat{g}_{\mathrm{can}}(\cdot,w)]$ 
is a Green's current satisfying the differential equation
\begin{align}\label{gcancurrent}
d_{z}d^{c}_{z}[\widehat{g}_{\mathrm{can}}(z,w)](f)+ f(w)=[\hatcan(z)](f),
\end{align}
where $f$ is a  smooth function defined on $\overline{X}$. We refer the reader to Theorem \Rmnum{2}.1.5 in \cite{lang} for the proof of the 
above equation. Let us denote the restriction of $\widehat{g}_{\mathrm{can}}(z,w)$ to $X\times X$ by $\gcan(z,w)$. 
\paragraph{Residual canonical metric on $\Omega_{\overline{X}}^{1}$}
Let $\Omega_{\overline{X}}^{1}$ denote the cotangent bundle of holomorphic differential forms on 
$\overline{X}$. For $z \in \overline{X}$, we define
\begin{equation*}
\|d\vartheta_{z}\|_{\mathrm{res,can}}^{2}(z)=\exp\bigg(\lim_{w\rightarrow z}
\big(\widehat{g}_{\mathrm{can}}(z,w)+ \log |\vartheta_{z}(w)|^{2}\big)\bigg).
\end{equation*}
From equation \eqref{gcanbounded}, it follows that the residual canonical metric is well defined and remains 
smooth on $\overline{X}.$ Furthermore, for $z \in \overline{X}$, the first Chern form $c_{1}\big(\Omega_{\overline{X}}^{1},\|\cdot\|_{\mathrm{res,can}}\big)$ is given by the formula
\begin{align}\label{chernformcan}
c_{1}\big(\Omega_{\overline{X}}^{1},\|\cdot\|_{\mathrm{res,can}}\big)= 
-d_{z}d_{z}^{c}\log\|d\vartheta_{z}\|_{\mathrm{res,can}}^{2}(z) =
(2g-2)\hatcan(z).
\end{align}
We refer the reader to \cite{arakelov} for the details of the proof of the above formula.
\paragraph{Parabolic Eisenstein Series} 
For $z\in X$ and $s\in\mathbb{C}$ with $\Re(s)> 1$, the parabolic Eisenstein series $\mathcal{E}_{\mathrm{par},p}(z,s)$ corresponding to a 
cusp $p\in\mathcal{P}$ is defined by the series
\begin{equation*}
\mathcal{E}_{\mathrm{par},p}(z,s) = \sum_{\eta \in \Gamma_{p}\backslash \Gamma}\Im(\sigma_{p}^{-1}\eta z)^{s}.
\end{equation*}
The series converges absolutely  and uniformly for $\Re(s) >1 $. It admits a meromorphic continuation to all 
$s\in\mathbb{C}$ with a simple pole at $s = 1$, and the Laurent expansion at $s=1$ is of the form 
\begin{equation}\label{laurenteisenpar}
\mathcal{E}_{\mathrm{par},p}(z,s) = \frac{1}{\vx(X)}\cdot\frac{1}{s-1} + 
\kappa_{p}(z) + O_{z}(s-1),
\end{equation}
where $\kappa_{p}(z)$ the constant term of $\mathcal{E}_{\mathrm{par},p}(z,s)$ at $s=1$ is called Kronecker's limit function (see Chapter 6 of 
\cite{hi}).
\paragraph{Heat Kernels}
For $t \in \mathbb{R}_{> 0}$  and $z, w \in \mathbb{H}$, the hyperbolic heat kernel $K_{\mathbb{H}}(t;z,w)$ on $\mathbb{R}_{> 0}\times
\mathbb{H} \times\mathbb{H}$ is given by the formula
\begin{equation*}\label{defnkh}
K_{\mathbb{H}}(t;z,w)= \frac{\sqrt{2}e^{- t\slash 4}}{(4\pi t)^{3\slash 2}}
\int_{d_{\mathbb{H}}(z,w)}^{\infty}\frac{re^{-r^{2}\slash 4t}}{\sqrt{\cosh(r)-\cosh (d_{\mathbb{H}}(z,w))}}dr,
\end{equation*}
where $d_{\mathbb{H}}(z,w)$ is the hyperbolic distance between $z$ and $w$.

For $t \in  \mathbb{R}_{> 0}$ and $z, w \in X$, the hyperbolic heat kernel $\khyp(t;z,w)$ on $\mathbb{R}_{> 0}\times X\times X$ is defined as 
\begin{equation*}\label{defnkhyp}
\khyp(t;z,w)=\sum_{\gamma\in\Gamma}K_{\mathbb{H}}(t;z,\gamma w).
\end{equation*}
For $z,w\in X,$ the hyperbolic heat kernel $\khyp(t;z,w)$ satisfies the differential equation
\begin{align}\label{diffeqnkhyp} 
\bigg(\Delta_{\text{hyp},z} + \frac{\partial}{\partial t}\bigg)\khyp(t;z,w) &=0,
\end{align}
Furthermore for a fixed $w\in X$ and any smooth function $f$ on $X$, the hyperbolic heat kernel $\khyp(t;z,w)$ 
satisfies the equation
\begin{align}\label{normcondkhyp}  
\lim_{t\rightarrow 0}\int_{X}\khyp(t;z,w)f(z)\hyp(z) &= f(w).
\end{align}
To simplify notation, we write $\khyp(t;z)$ instead of $\khyp(t;z,z)$, when $z=w$.
\paragraph{Automorphic Green's function}
For $z, w \in \mathbb{H}$ with $z \not= w$, and  $s$ $\in$ $\mathbb{C}$ with $\Re(s)> 0$, the free-space Green's function $g_{\mathbb{H},s}(z,w)$ is defined as
\begin{equation*} 
g_{\mathbb{H},s}(z,w) = g_{\mathbb{H},s}(u(z,w))= \dfrac{\Gamma(s)^{2}}{\Gamma(2s)}u^{-s}
F(s,s;2s,-1\slash u),
\end{equation*} 
where $u=u(z,w)=|z-w|^{2}\slash( 4\Im(z)\Im(w))$ and $F(s,s;2s,-1\slash u)$ is the 
hypergeometric function. 

There is a sign error in the formula defining the free-space Green's function given by equation (1.46) in \cite{hi}, i.e., 
the last argument $-1\slash u$ in the hypergeometric function has been incorrectly stated as $1\slash u$, which we have 
corrected in our definition. We have also normalized the free-space Green's function defined in \cite{hi} by multiplying it by 
$4\pi.$

\vspace{0.2cm}
For $z, w \in X$ with $z\not = w$, and $s\in\mathbb{C}$ with $\Re(s) > 1$, the automorphic Green's function 
$g_{\mathrm{hyp},s}(z,w)$ is defined as
\begin{equation*}
g_{\mathrm{hyp},s}(z,w) = \sum_{\gamma\in\Gamma}g_{\mathbb{H},s}(z,\gamma w).
\end{equation*}
The series converges absolutely uniformly for $z\not = w$ and $\Re(s) > 1$ (see Chapter 5 in \cite{hi}). 

For $z, w \in X$ with $z \not = w$, and $s\in\mathbb{C}$ with $\Re(s) > 1$, the automorphic Green's function satisfies the 
following properties (see Chapters 5 and 6 in \cite{hi}):

(1) For $\Re(s(s-1)) > 1$, we have
\begin{equation*}
g_{\mathrm{hyp},s}(z,w) = 4\pi\int_{0}^{\infty}\khyp(t;z,w)e^{-s(s-1)t}dt.
\end{equation*}
(2) It admits a logarithmic singularity along the diagonal, i.e.,  
\begin{equation*}
\lim_{w\rightarrow z}\big(g_{\mathrm{hyp},s}(z,w) + \log{|\vartheta_{z}(w)|^{2}}\big)= O_{s,z}(1).
\end{equation*}
(3) The automorphic Green's function $g_{\mathrm{hyp},s}(z,w)$ admits a meromorphic continuation to all $s\in\mathbb{C}$ with a 
simple pole at $s=1$ with residue $4\pi\slash\vx(X)$, and the Laurent expansion at $s=1$ is of the form
\begin{equation*}
g_{\mathrm{hyp},s}(z,w)= \frac{4\pi}{s(s-1)\vx(X)} + g^{(1)}_{\mathrm{hyp}}(z,w) + O_{z,w}(s-1),
\end{equation*} 
where $g_{\mathrm{hyp}}^{(1)}(z,w)$ is the constant term of $g_{\mathrm{hyp},s}(z,w)$ at $s=1$. 

\vspace{0.15cm}
(4) Let $p,q\in \mathcal{P}$ be two cusps. Put
 \begin{align*}
C_{p,q}  = \min \bigg{\lbrace} c > 0\,\bigg{|} \bigg(\begin{array}{ccc} a &b\\
 c & d  \end{array}\bigg) \in \sigma_{p}^{-1}\Gamma \sigma_{q} \bigg{\rbrace},
\end{align*} 
and $C_{p,p}=C_{p}$. Then, for $z,w\in X$ with $\Im(w) > \Im(z)$ and  $\Im(w)\Im(z) > C_{p,q}^{-2}$, and $s\in\mathbb{C}$ with $\Re(s) > 1$, 
the automorphic Green's function admits the Fourier expansion
\begin{align} 
g_{\mathrm{hyp},s}(\sigma_{p}z,\sigma_{q}w)= \frac{4\pi\Im(w)^{1-s}}{2s-1}
\mathcal{E}_{\mathrm{par},q}(\sigma_{p}z,s) -\delta_{p,q} \log\big|1- e^{2\pi i (w-z)}\big|^{2} 
+O\big(e^{-2\pi (\Im( w)-\Im (z))}\big). \label{fourierautghyp}
\end{align}
This equation has been proved as Lemma 5.4 in \cite{hi}, and one of the terms was wrongly estimated in the proof of the lemma. 
We have corrected this error, and stated the corrected equation. 
\paragraph{Hyperbolic Green's function}
For $z, w \in X$ and $z\not = w$, the hyperbolic Green's function is defined as 
\begin{equation*}
\ghyp(z,w) = 4\pi\int_{0}^{\infty}\bigg(\khyp(t;z,w)-\frac{1}{\vx(X)}\bigg)dt.
\end{equation*}
For $z, w \in X$ with $z \not = w$, the hyperbolic Green's function satisfies the 
following properties:

(1) For $z, w \in X$, we have 
\begin{equation}\label{ghypbounded}
\lim_{w\rightarrow z}\big( \ghyp(z,w) + \log{|\vartheta_{z}(w)|^{2}}\big)= O_{z}(1).
\end{equation}
(2) For $z, w \in X\backslash \mathcal{E}$, the hyperbolic Green's function satisfies 
the differential equation 
\begin{align}
d_{z}d_{z}^{c}\ghyp(z,w) +\delta_{w}(z)& = \shyp(z), \label{diffeqnghyp} \\
\intertext{ with the normalization condition} 
\int_{X}\ghyp(z,w)\hyp(z) & = 0. \label{normcondghyp}
\end{align}
(3) For $z,w\in X$ and $z\not=w$, we have
\begin{equation}\label{laurentghyp}
 \ghyp(z,w)= g^{(1)}_{\mathrm{hyp}}(z,w)= \lim_{s\rightarrow 1}\bigg(g_{\text{hyp},s}(z,w) - 
\frac{4\pi}{s(s-1)\vx(X)}\bigg).
\end{equation}
The above properties follow from the properties of the heat kernel $\khyp(t;z,w)$ (equations \eqref{diffeqnkhyp} and \eqref{normcondkhyp}) or 
from that of the automorphic Green's function $g_{\mathrm{hyp},s}(z,w)$. 
\paragraph{Residual hyperbolic metric on $\Omega_{X}^{1}$}
For $z\in X$, we define 
\begin{equation*}
\|d\vartheta_{z}\|^{2}_{\mathrm{res,hyp}}(z)=\exp\bigg(\lim_{w\rightarrow z}
\big(\ghyp(z,w)+\log|\vartheta_{z}(w)|^{2}\big)\bigg).
\end{equation*} 
From equation \eqref{ghypbounded}, it follows that the residual hyperbolic metric is well defined on $X$. Furthermore, 
from Proposition 3.3 in \cite{jk}, for $z\in X\backslash \mathcal{E}$, we have  
\begin{align}
c_{1}\big(\Omega_{X}^{1},\|\cdot\|_{\mathrm{res,hyp}}\big)=&\,-d_{z}d_{z}^{c}\|
d\vartheta_{z}\|^{2}_{\mathrm{res,hyp}}(z)=\frac{1}{2\pi}\hyp(z)+ \bigg(\int_{0}^{\infty}
\del \khyp(t;z)dt\bigg)\hyp(z).\label{chernformhyp}
\end{align}
\paragraph{Convergence results}
From Lemmas 5.2 and 6.3, and Proposition 7.3 in \cite{K}, the function
\begin{align*}
4\pi\int_{0}^{\infty}\del \khyp(t;z)dt
\end{align*}
is well defined on $X$ and remains bounded at the elliptic fixed points and at the cusps. Hence, it defines a 
smooth function on $\overline{X}$, which we denote symbolically by
\begin{equation*}
\int_{0}^{\infty}\del \widehat{K}_{\mathrm{hyp}}(t;z)dt.
\end{equation*}
\paragraph{Key identity}
For $z \in X\backslash \mathcal{E}$, we have the relation of differential forms
\begin{align}\label{keyidentity}
&g\can(z) =\bigg(\frac{1}{4\pi}+ \frac{1}{\vx(X)}\bigg)\hyp(z)+ \frac{1}{2}\bigg(\int_{0}^{\infty}\del \khyp(t;z) dt\bigg)\hyp(z).
\end{align}
This relation has been established as Theorem 3.4 in \cite{jk}, when $X$ is 
compact. The proof given in \cite{jk} applies to our case where $X$ does admit elliptic fixed points and 
cusps, as long as $z\in X\backslash \mathcal{E}.$
\paragraph{The space $C_{\ell,\ell\ell}(\overline{X})$}
Let $C_{\ell,\ell\ell}(\overline{X})$ denote the set of complex-valued functions $f:\overline{X}\rightarrow 
\mathbb{P}^{1}(\mathbb{C})$, which admit the following type of singularities at finitely many points 
$\mathrm{Sing}(f)\subseteq \overline{X}$, and are smooth away from $\mathrm{Sing}(f)$: 

(1) If $s\in\mathrm{Sing}(f)\backslash \mathcal{P}$, then as $z$ approaches $s$, the function $f$ satisfies  
\begin{align}\label{fsingular}
f(z)= c_{f,s}\log|\vartheta_{s}(z)| + O_{z}(1),
\end{align}
for some $c_{f,s}\in \mathbb{C}$.

(2) For $p\in\mathrm{Sing}(f)\cap 
\mathcal{P}$, as $z$ approaches $p$, the function $f$ satisfies
\begin{align}\label{fcusp}
f(z)=c_{f,p}\log\big(-\log|\vartheta_{p}(z)|\big) + O_{z}(1),
\end{align}
for some $c_{f,p}\in \mathbb{C}$. 
\section{Extension of key identity}\label{section2}
In this section, we extend the key identity, i.e., equation \eqref{keyidentity} to elliptic fixed points and 
cusps at the level of currents acting on the space of singular functions $C_{\ell,\ell\ell}(\overline{X})$. 

In subsection \ref{subsection2.1}, we investigate the behavior of the hyperbolic Green's function at the cusps, 
and show that it defines a current acting on the space of singular functions $C_{\ell,\ell\ell}(\overline{X})$. In subsection 
\ref{subsection2.2}, we prove an auxiliary identity which is useful in extending the key identity 
\eqref{keyidentity} to elliptic fixed points and cusps. In subsection \ref{subsection2.3}, using the results 
from the previous two subsections, we extend the key identity.  
\subsection{Hyperbolic Green's function as a Green's current}\label{subsection2.1}
Although it is obvious from the differential equation \eqref{diffeqnghyp} that $\ghyp(z,w)$ is 
$\log\log$-singular at the cusps, the exact asymptotics derived in the following proposition come very useful 
in the upcoming articles (especially in \cite{anil}). 
\begin{prop}\label{prop1}
With notation as in Section \ref{section1}, for a fixed $ w\in X$, and for $z\in X$ with 
$\Im(\sigma_{p}^{-1}z)>\Im(\sigma_{p}^{-1} w)$ and $\Im(\sigma_{p}^{-1}z)\Im(\sigma_{p}^{-1}w) >C_{p}^{-2}$, 
we have
\begin{align}
& \ghyp(z,w) = 4\pi\kappa_{p}(w) - 
\frac{4\pi}{\vx(X)}-\frac{4\pi\log\big(\Im(\sigma_{p}^{-1}z)\big)}{\vx(X)}-\notag \\ &
 \log\big{|}1-e^{2\pi i(\sigma_{p}^{-1}z - \sigma_{p}^{-1}w)}\big{|}^{2}+
O\big(e^{-2\pi (\Im(\sigma_{p}^{-1}z)-\Im(\sigma_{p}^{-1}w))}\big).
\label{propeqn}
\end{align}
\begin{proof}
As the limit in \eqref{laurentghyp} converges uniformly, combining it with equation \eqref{fourierautghyp}, for a fixed $w \in X$, 
for each $z$ $\in$ $X$ with $\Im(\sigma_{p}^{-1}z) > \Im(\sigma_{p}^{-1}w)$ 
and $\Im(\sigma_{p}^{-1}z)\Im(\sigma_{p}^{-1}w)>C_{p}^{-2}$, and $s\in\mathbb{C}$ with $\Re(s) > 1$, we have
\begin{align}
\ghyp(z,w)\,= \,&4\pi\lim_{s\rightarrow 1}\bigg(\frac{\Im(\sigma_{p}^{-1}z)^{1-s}}{2s-1}\mathcal{E}_{\mathrm{
par},p}(w,s)-\frac{1}{(s-1)\vx(X)}\bigg)+\frac{4\pi}{\vx(X)} - \notag\\&\log\big{|}1 - 
e^{2\pi i(\sigma_{p}^{-1}z- \sigma_{p}^{-1}w)}\big{|}^{2} + O\big(e^{-2\pi(\Im(\sigma_{p}^{-1}z)-
\Im(\sigma_{p}^{-1}w))}\big).
\label{proofhypcuspgreen1}
\end{align}
To evaluate the above limit, we compute the Laurent expansions of 
$\mathcal{E}_{\mathrm{par},p}(w,s)$, $\Im(\sigma_{p}^{-1}z)^{1-s}$, and 
$(2s-1)^{-1}$ at $s=1$. The Laurent expansions 
of $\Im{(\sigma_{p}^{-1}z)^{1-s}}$ and $(2s-1)^{-1}$ at $s=1$ are easy to compute, 
and are of the form
\begin{align*}
&\Im{(\sigma_{p}^{-1}z)}^{1-s}=1 - (s-1)\log\big(\Im{(\sigma_{p}^{-1}z)}\big)
+ O\big((s-1)^{2}\big), \,\,\,\,\frac{1}{2s-1}  = 1- 2(s-1) + O\big((s-1)^{2}\big).
\end{align*}
Combining the above two equations with equation \eqref{laurenteisenpar}, we arrive at 
\begin{align*}
&4\pi\lim_{s\rightarrow 1}\bigg(\frac{\Im(\sigma_{p}^{-1}z)^{1-s}}{2s-1}
\mathcal{E}_{\mathrm{par},p}(w,s)-\frac{1}{(s-1)\vx(X)}\bigg)=\notag\\&
 4\pi \kappa_{p}(w) -\frac{8\pi}{\vx(X)}- \frac{4\pi\log\big(\Im(\sigma_{p}^{-1}z)\big)}{\vx(X)},
\end{align*}
which together with equation (\ref{proofhypcuspgreen1}) implies the proposition.
\end{proof}
\end{prop}  
\begin{cor}\label{cor2}
For a fixed $w\in X$, as $z\in X$ approaches a cusp $p\in\mathcal{P}$, we have
\begin{align*}
\ghyp(z,w) & = -\frac{4\pi\log\big(\Im(\sigma_{p}^{-1}z)\big)}{\vx(X)}+
O_{z,w}(1) = -\frac{4\pi\log\big(-\log|\vartheta_{p}(z)|\big)}{\vx(X)}+ O_{z,w}(1).
\end{align*}
\begin{proof}
The corollary follows from Proposition \ref{prop1}.
\end{proof}
\end{cor}
From the above corollary, it follows that for a fixed $w\in X$, as a function in the variable $z$, the 
hyperbolic Green's function $\ghyp(z,w)$ has $\log\log$-growth at the cusps. Hence, for a fixed 
$w\in\overline{X}\backslash \mathcal{P}$, as a function in the variable $z$, it defines a singular function 
$\widehat{g}_{\mathrm{hyp}}(z,w)$ on $\overline{X}$ with $\log\log$-singularity cusps and $\log$-singularity 
at $z=w$. So for a fixed $w\in\overline{X}$,  the hyperbolic Green's function $g_{\mathrm{hyp}}(z,w)$ 
determines a current $[\widehat{g}_{\mathrm{hyp}}(\cdot,w)]$ of type (0,0) acting on smooth (1,1)-forms defined on $\overline{X}$. 
\begin{rem}
For any $f\in C_{\ell,\ell\ell}(\overline{X})$, from standard arguments from analysis, it follows 
that $d_{z}d_{z}^{c}f(z)$ defines an integrable (1,1)-form on $\overline{X}$. Furthermore, for a fixed 
$w\in\overline{X}\backslash(\mathrm{Sing}(f)\cup \mathcal{P})$, the following integral exists 
\begin{align*}
\int_{\overline{X}}\widehat{g}_{\mathrm{hyp}}(z,w)d_{z}d_{z}^{c}f(z).
\end{align*}
\end{rem}
In the following lemma, we show that the hyperbolic Green's function defines a Green's current acting on the 
space of singular functions $C_{\ell,\ell\ell}(\overline{X})$.  
\begin{lem}\label{ghyplemma}\label{lem3}
Let $f\in C_{\ell,\ell\ell}(\overline{X})$, then for a $w\in\overline{X}\backslash(\mathrm{Sing}(f)\cup 
\mathcal{P})$ fixed, we have the equality of integrals
\begin{align*}
\int_{\overline{X}}\widehat{g}_{\mathrm{hyp}}(z,w)d_{z}d_{z}^{c}f(z) + f(w)+ 
\sum_{\substack{s\in \mathrm{Sing}(f)\\s\not\in \mathcal{P}}}
\frac{c_{f,s}}{2}\widehat{g}_{\mathrm{hyp}}(s,w)=
\int_{\overline{X}}f(z)\hatshyp(z).
\end{align*}
\begin{proof}
Let $w\in\overline{X}\backslash(\mathrm{Sing}(f)\cup \mathcal{P})$ and let $U_{r}(w)$, $U_{r}(s)$, and 
$U_{r}(p)$ denote open coordinate disks of radius $r$ around $w$, $s\in\mathrm{Sing}(f)$, and a cusp 
$p\in \mathcal{P}$, respectively. Put
\begin{equation*}
Y_{r} = \overline{X}\backslash \bigg(U_{r}(w)\cup\bigcup_{\substack{s\in \mathrm{Sing}(f)\\s\not\in
\mathcal{P}}}U_{r}(s)\cup\bigcup_{p\in \mathcal{P}}U_{r}(p)\bigg).
\end{equation*} 
From equation (\ref{diffeqnghyp}) and Stokes's theorem, it follows that it suffices to prove that 
\begin{align*}
&\int_{Y_{r}}\ghyp(z,w)d_{z}d_{z}^{c}f(z)-  \int_{Y_{r}}f(z)\shyp(z)
=\notag\\[0.3em]&
\int_{\partial  U_{r}(w)}\ghyp(z,w)(-d_{z}^{c}f(z)) -
\int_{\partial U_{r}(w)}f(z)(-d_{z}^{c}\ghyp(z,w))+\notag
\\&\sum_{\substack{s\in  \mathrm{Sing}(f)\\s\not\in \mathcal{P}}}
\bigg(\int_{\partial U_{r}(s)} \ghyp(z,w) (-d_{z}^{c}f(z))- 
\int_{\partial U_{r}(s)}f(z)(-d_{z}^{c}\ghyp(z,w))\bigg)
 +\notag\\&\hspace{0.4cm}\sum_{p\in \mathcal{P}}\bigg(
\int_{\partial U_{r}(p)}\ghyp(z,w)(-d_{z}^{c}f(z))-
\int_{\partial U_{r}(p)}f(z)(-d_{z}^{c}\ghyp(z,w))\bigg)\xrightarrow[r\rightarrow 0]{}\\ &-f(w)
-\frac{c_{f,s}}{2}\widehat{g}_{\mathrm{hyp}}(s,w).
\end{align*}
Recall that $d_{z}^{c}$ in polar coordinates is given by  
\begin{align*}
 d_{z}^{c}= \frac{r}{2}\frac{\partial}{\partial r}\frac{d\theta}{2\pi} - 
\frac{1}{{4\pi}}\frac{\partial}{\partial \theta}\frac{dr}{r}.
\end{align*}
Hence, as $w\not\in\mathrm{Sing}(f)$, using equation \eqref{ghypbounded} we derive 
\begin{align}
&\int_{\partial U_{r}(w)}\ghyp(z,w)(-d^{c}f(z)) -\int_{\partial U_{r}(w)}f(z)(-d_{z}^{c}\ghyp(z,w))=\notag\\
&\int_{0}^{2\pi}r\log r\frac{\partial f}{\partial r} \frac{d\theta}{2\pi}-\int_{0}^{2\pi}f(z)r\frac{\partial \log r}{\partial r}\frac{d\theta}{2\pi} + O(r)
\xrightarrow[r\rightarrow 0]{}-f(w).\label{lem3eqn2}
\end{align}
As $w\not\in  \mathrm{Sing}(f)$, the hyperbolic Green's function $\ghyp(z,w)$ remains smooth at 
$s\in  \mathrm{Sing}(f)\backslash\mathcal{P}$. So for any $s\in\mathrm{Sing}(f)$ and $s\not\in\mathcal{P}$, 
using equation (\ref{fsingular}) and from similar computations as in \eqref{lem3eqn2}, we get
\begin{align}
&\int_{\partial U_{r}(s)}\ghyp(z,w)(-d_{z}^{c}f(z))-\int_{0}^{2\pi}f(z)(-d_{z}^{c}\ghyp(z,w))=\notag\\[0.5em]
&-c_{f,s}\bigg(\int_{0}^{2\pi}\ghyp(z,w)\frac{r}{2}\frac{\partial \log r}{\partial r}
\frac{d\theta}{2\pi} +\int_{0}^{2\pi}\frac{r}{2}\log r\frac{\partial \ghyp(z,w)}{\partial r}\frac{d\theta}{2\pi}\bigg) + O(r)
\xrightarrow[r\rightarrow 0]{}-\frac{c_{f,s}}{2}\ghyp(s,w).\label{lem3eqn3}
\end{align}
Finally for any cusp $p\in\mathcal{P}$, using Corollary \ref{cor2} and equation \eqref{fcusp}, we compute
\begin{align}
&\int_{\partial U_{r}(p)}\ghyp(z,w)(-d_{z}^{c}f(z))-\int_{\partial U_{r}(p)}f(z)(-d_{z}^{c}\ghyp(z,w)) = 
\notag \\&\frac{4\pi c_{f,p}}{\vx(X)}\bigg(\int_{0}^{2\pi}\log\big{(}-\log r\big{)}\frac{r}{2}\frac{\partial
\log\big{(}-\log r\big{)}}{\partial r}\frac{d\theta}{2\pi} -\int_{0}^{2\pi}\log
\big{(}-\log r\big{)}\frac{r}{2}\frac{\partial\log\big{(}-\log r\big{)}}{\partial r}\frac{d\theta}{2\pi}\bigg) +
\notag \\&O(1\slash \log r) = O(1\slash \log r)\xrightarrow[r\rightarrow 0]{} 0.\label{lem3eqn4}
\end{align}
Combining equations \eqref{lem3eqn2}, \eqref{lem3eqn3}, and \eqref{lem3eqn4} completes the proof of the lemma. 
\end{proof}
\end{lem}
\begin{cor}\label{cor4}
Let $f\in C_{\ell,\ell\ell}(\overline{X})$, then for a fixed $w\in X\backslash \left(\mathrm{Sing}(f)\cap X\right)$, 
we have the equality of integrals
\begin{align*}
\int_{X}\ghyp(z,w)d_{z}d_{z}^{c}f(z) + f(w)+ 
\sum_{\substack{s\in \mathrm{Sing}(f)\\s\not\in \mathcal{P}}}
\frac{c_{f,s}}{2}\ghyp(s,w)= \int_{X}f(z)\shyp(z).
\end{align*}
\begin{proof}
The proof follows from Lemma \ref{lem3} and the fact that there are only finitely many cusps of $X$.   
\end{proof}
\end{cor}
\subsection{An auxiliary identity}\label{subsection2.2}
In this subsection, we drive an auxiliary identity, which is useful in proving the key identity in next 
subsection.
\begin{lem}\label{lem5}
There exists a unique integrable function $\widehat{\Phi}(z)$ defined on $\overline{X}$, which satisfies the 
differential equation
\begin{align}
& d_{z}d_{z}^{c}[\widehat{\Phi}(z)]=[\hatshyp(z)]-[\hatcan(z)], \label{lem5eqn1}
\intertext{with the normalization condition} 
& \int_{X}\widehat{\Phi}(z)\can(z) =0, \label{lem5eqn2} 
\end{align}
where $[\widehat{\Phi}(z)]$ is the current determined by $\widehat{\Phi}(z)$ acting on smooth (1,1)-forms 
defined on $\overline{X}$.
\begin{proof}
Since the cohomology classes of $[\hatshyp(z)]$ and $[\hatcan(z)]$ are equal in  $H^{2}\left(\overline{X},\mathbb{Z}\right)\cong \mathbb{Z}$, the difference 
$[\hatshyp(z)]-[\hatcan(z)]$ is a $d$-exact current on $\overline{X}$. Hence, from the $\partial\bar{\partial}$-lemma for currents, we can conclude that there 
exists an integrable function $\widehat{\Phi}(z)$ defined on $\overline{X}$ such that
\begin{equation*}
d_{z}d_{z}^{c}[\widehat{\Phi}(z)]= [\hatshyp(z)]-[\hatcan(z)],
\end{equation*}
which proves the existence of $\widehat{\Phi}(z)$. The normalization condition \eqref{lem5eqn2} ensures the uniqueness of $\widehat{\Phi}(z)$.
\end{proof}
\end{lem}
\begin{lem}\label{lem6}
Let us denote the restriction of $\widehat{\Phi}(z)$ to $X$ by $\Phi(z)$. Then, for $z, w\in X$, we have
\begin{align}\label{lem6eqn}
\ghyp(z,w)-\gcan(z,w)= \frac{1}{2}\bigg(\Phi(z) + \Phi(w) + \int_{X}\ghyp(z,\zeta)\can(\zeta)+ 
\int_{X}\ghyp(w,\zeta)\can(\zeta)\bigg).
\end{align}
\begin{proof}
For a fixed $w \in X$, consider the function 
\begin{equation*}
F_{w}(z)= \ghyp(z,w)-\gcan(z,w)-\int_{X}\ghyp(w,\zeta)\can(\zeta)
\end{equation*}
defined on $X$. As $\gcan(z,w)$ and $\ghyp(z,w)$ define currents $[\widehat{g}_{\mathrm{can}}(\cdot,w)]$ and 
$[\widehat{g}_{\mathrm{hyp}}(\cdot,w)]$ of type (0,0) on $\overline{X}$, respectively, the function 
$F_{w}(z)$ determines a current  
\begin{equation*}
[\widehat{F}_{w}]= [\widehat{g}_{\mathrm{hyp}}(\cdot,w)]-[\widehat{g}_{\mathrm{can}}(\cdot,w)] -\int_{X}\ghyp(w,\zeta)\can(\zeta)
\end{equation*}
of type (0,0) acting on smooth smooth (1,1)-forms defined on $\overline{X}$. For a fixed $w\in X$, using 
equation \eqref{gcancurrent} and Lemma \ref{lem3}, it is easy to see that $\widehat{F}_{w}$ 
satisfies equations \eqref{lem5eqn1} and \eqref{lem5eqn2}. Hence, from the uniqueness of 
$\widehat{\Phi}(z)$, we get
\begin{equation*}
\Phi(z) = F_{w}(z)= \ghyp(z,w)-\gcan(z,w)-\int_{X}\ghyp(w,\zeta)\can(\zeta),
\end{equation*}
which implies that $F_{w}(z)$ is independent of $w\in X$. Hence, from the above equation and from the symmetry of the Green's functions 
$\ghyp(z,w)$ and $\gcan(z,w)$, we deduce that
\begin{align*}
\ghyp(z,w)-\gcan(z,w)= \frac{1}{2}\bigg(\Phi(z) + \Phi(w)+  \int_{X}\ghyp(z,\zeta)\can(\zeta)+ 
\int_{X}\ghyp(w,\zeta)\can(\zeta)\bigg),
\end{align*}
which proves the lemma.
\end{proof}
\end{lem}
\begin{prop}\label{prop7}
For $z,w\in X$, we have
\begin{equation*}
\ghyp(z,w)-\gcan(z,w)= \phi(z) + \phi(w),
\end{equation*}
where
\begin{align*} 
\phi(z)& = \int_{X}\ghyp(z,\zeta)\can(\zeta) -\frac{1}{2} \int_{X}\int_{X} \ghyp(\xi,\zeta)\can(\zeta)\can(\xi).
\end{align*}
\begin{proof}
For all $z,w\in X$, combining Lemma \ref{lem6} and equations \eqref{normcondgcan} and \eqref{lem5eqn2}, we obtain
\begin{align*}
2&\int_{X}\big(\ghyp(z,w)-\gcan(z,w)\big)\can(w) =2\int_{X}\ghyp(z,w)\can(w)=\notag\\&\int_{X}\bigg(\Phi(z) +
\Phi(w) + \int_{X}\ghyp(z,\zeta)\can(\zeta)+ \int_{X}\ghyp(w,\zeta)\can(\zeta)\bigg)\can(w)=\notag\\&
\Phi(z)+\int_{X}\ghyp(z,\zeta)\can(\zeta)+\int_{X}\int_{X}\ghyp(\xi,\zeta)\can(\zeta)\can(\xi).
\end{align*}
Hence, we arrive at
\begin{align*}
\Phi(z)= \int_{X}\ghyp(z,\zeta)\can(\zeta)-\int_{X}\int_{X}\ghyp(\xi,\zeta)\can(\zeta)\can(\xi). 
\end{align*}
Substituting the above formula for $\Phi(z)$ in equation (\ref{lem6eqn}), we get 
\begin{align*}
& \ghyp(z,w)-\gcan(z,w)= \\ 
&\int_{X}\ghyp(z,\zeta)\can(\zeta)+\int_{X}\ghyp(w,\zeta)\can(\zeta) -
\int_{X}\int_{X}\ghyp(\xi,\zeta)\can(\zeta)\can(\xi).
\end{align*}
The proof of the proposition follows by setting
\begin{align*}
\phi(z) =\int_{X}\ghyp(z,\zeta)\can(\zeta)-\frac{1}{2}
\int_{X}\int_{X}\ghyp(\xi,\zeta)\can(\zeta)\can(\xi).
\end{align*}  
\end{proof}
\end{prop}
\begin{cor}\label{cor8}
As $z\in X$ approaches a cusp $p\in \mathcal{P}$, we have 
\begin{equation*}
\phi(z)=  -\frac{4\pi\log\big(-\log|\vartheta_{p}(z)|\big)}{\vx(X)}+O_{z}(1).
\end{equation*}
\begin{proof}
For a fixed $w\in X$, as a function in the variable $z$, the canonical Green's function 
$\gcan(z,w)$ remains smooth at the cusps. So the proof of the corollary follows directly from combining 
Proposition \ref{prop7} and Corollary \ref{cor2}.  
\end{proof}
\end{cor}
In the following proposition, we show that the residual hyperbolic metric is $\log\log$-singular at the cusps.
\begin{cor}\label{cor9}
As $z\in X$ approaches a cusp $p\in \mathcal{P}$, we have
\begin{align*}
\log\|d\vartheta_{z}\|^{2}_{\mathrm{res,hyp}}(z)=-\frac{8\pi\log\big(-\log|\vartheta_{p}(z)|\big)}{\vx(X)} + 
O_{z}(1).
\end{align*}
\begin{proof}
From Proposition \ref{prop7}, we have 
\begin{align*}
\lim_{w\rightarrow z}\big(\ghyp(z,w)+ \log |\vartheta_{z}(w)|^{2}\big)= \lim_{w\rightarrow z}
\big(\gcan(z,w)+ \log |\vartheta_{z}(w)|^{2}\big) + 2\phi(z).
\end{align*}
The proof of the corollary follows directly from combining equation \eqref{gcanbounded} and Corollary \ref{cor8}.
\end{proof}
\end{cor}
From Corollary \ref{cor9}, it follows that $\log\|d\vartheta_{z}\|^{2}_{\mathrm{res,hyp}}(z)$ is smooth on $\overline{X}\backslash
\mathcal{P}$, and admits a $\log\log$-singularity at the cusps. So for any 
$f\in C_{\ell,\ell\ell}(\overline{X})$, the following integral exists
\begin{align*}
\int_{\overline{X}}\log\|d\vartheta_{z}\|_{\mathrm{res,hyp}}^{2}(z)d_{z}d_{z}^{c}f(z) .
\end{align*}
\subsection{Key identity}\label{subsection2.3}
In this subsection, we extend relation \eqref{keyidentity} to elliptic fixed points and cusps at the level of 
currents acting on the space of singular functions $C_{\ell,\ell\ell}(\overline{X})$. 

Let $U_{r_{0}}(s)$ denote an open coordinate disk of fixed radius $r_{0}$ 
around $s\in \mathrm{Sing}(f)\cup \mathcal{S}$, and $r_{0}$ is small enough such that any two coordinate disks are disjoint. Put
\begin{align*}
U_{r_{0}}=\bigcup_{s\in  \mathrm{Sing}(f)\cup \mathcal{S}}U_{r_{0}}(s)\,\,\,\mathrm{and}\,\,\, Y_{r_{0}}= \overline{X}\backslash  U_{r_{0}}.
\end{align*}
Furthermore, for $0< r< r_{0}$, let $U_{r}(s)$ denote an open coordinate disk of radius $r$ around $s\in \mathrm{Sing}(f)\cup\mathcal{S}$, and let 
$U_{r_{0},r}(s)$ denote the annulus $U_{r_{0}}(s)\backslash U_{r}(s)$. Put 
\begin{align*}
U_{r}=\bigcup_{s\in\mathrm{Sing}(f)\cup\mathcal{S}}U_{r}(s)\,\,\,
\mathrm{and}\,\,\,U_{r_{0},r}=U_{r_{0}}\backslash U_{r}.
\end{align*}
\begin{prop}\label{2.3prop1}
Let $f\in C_{\ell,\ell \ell}(\overline{X})$, then we have the equality of integrals
\begin{align}
-&\int_{U_{r_{0}}}\log\|d\vartheta_{z}\|^{2}_{\mathrm{res,can}}(z)
d_{z}d_{z}^{c}f(z) =(2g-2)\int_{U_{r_{0}}}f(z)\hatcan(z) + \sum_{\substack{s\in \mathrm{Sing}(f)\\s\not\in \mathcal{P}}}
\frac{c_{f,s}}{2}\log\|d\vartheta_{z}\|^{2}_{\mathrm{res,can}}(s)-\notag\\&
\int_{\partial U_{r_{0}}}\log\|d\vartheta_{z}\|^{2}_{\mathrm{res,can}}(z)
d_{z}^{c}f(z)+\int_{\partial U_{r_{0}}}f(z)d_{z}^{c}\log\|
d\vartheta_{z}\|^{2}_{\mathrm{res,can}}(z)\label{2.3prop1eqn}
\end{align}
\begin{proof}
From equation \eqref{chernformcan}, it follows that for any $r>0$, it suffices to prove that
\begin{align*}
-&\int_{U_{r_{0},r}}\log\|d\vartheta_{z}\|^{2}_{\mathrm{res,can}}
(z)d_{z}d_{z}^{c}f(z)+\int_{U_{r_{0},r}}f(z)d_{z}d_{z}^{c}\log\|d\vartheta_{z}\|^{2}_{\mathrm{res,can}}(z)
\xrightarrow[r\rightarrow 0]{}  \\&  \sum_{\substack{s\in \mathrm{Sing}(f)\\s\not\in \mathcal{P}}}
\frac{c_{f,s}}{2}\log\|d\vartheta_{z}\|^{2}_{\mathrm{res,can}}(s)-
\notag\\&\int_{\partial U_{r_{0}}}\log\|d\vartheta_{z}\|^{2}_{\mathrm{res,can}}(z)d_{z}^{c}f(z)+
\int_{\partial U_{r_{0}}}f(z)d_{z}^{c}\log\|d\vartheta_{z}\|^{2}_{\mathrm{res,can}}(z).
\end{align*} 
Using Stokes's theorem, we find that the left-hand side of the above limit simplifies to the following 
expression
\begin{align*}
&\int_{\partial U_{r}}\log\|d\vartheta_{z}\|^{2}_{\mathrm{res,can}}(z)d_{z}^{c}f(z)-
\int_{\partial U_{r}}f(z)d_{z}^{c}\log\|d\vartheta_{z}\|^{2}_{\mathrm{res,can}}(z)-\\
&\int_{\partial U_{r_{0}}}\log\|d\vartheta_{z}\|^{2}_{\mathrm{res,can}}(z)d_{z}^{c}f(z)+
\int_{\partial U_{r_{0}}}f(z) d_{z}^{c}\log\|d\vartheta_{z}\|^{2}_{\mathrm{res,can}}(z)
\end{align*} 
Furthermore, from the construction of the open sets $U_{r}$, we have
\begin{align*}
&\int_{\partial U_{r}}\log\|d\vartheta_{z}\|^{2}_{\mathrm{res,can}}(z)d_{z}^{c}f(z)-
\int_{\partial U_{r}}f(z)d_{z}^{c}\log\|d\vartheta_{z}\|^{2}_{\mathrm{res,can}}(z)=\\
&\sum_{\substack{s\in \mathrm{Sing}(f)\\s\not\in \mathcal{P}}}\bigg(\int_{\partial U_{r}(s)}\log\|d\vartheta_{z}\|^{2}_{\mathrm{res,can}}(z)d_{z}^{c}f(z)-
\int_{\partial U_{r}(s)}f(z)d_{z}^{c}\log\|d\vartheta_{z}\|^{2}_{\mathrm{res,can}}(z)\bigg)+\notag\\
&\,\,\,\,\,\,\sum_{p\in\mathcal{P}}\bigg(\int_{\partial U_{r}(p)}\log\|d\vartheta_{z}\|^{2}_{\mathrm{res,can}}(z)d_{z}^{c}f(z)-
\int_{\partial U_{r}(p)}f(z)d_{z}^{c}\log\|d\vartheta_{z}\|^{2}_{\mathrm{res,can}}(z)\bigg).
\end{align*}
As $\log\|d\vartheta_{z}\|^{2}_{\mathrm{res,can}}(z)$ remains smooth on $\overline{X}$, from arguments as in Lemma \ref{lem3}, we derive
\begin{align*}
&\sum_{\substack{s\in \mathrm{Sing}(f)\\s\not\in \mathcal{P}}}\bigg(\int_{\partial U_{r}(s)}\log\|d\vartheta_{z}\|^{2}_{\mathrm{res,can}}(z)d_{z}^{c}f(z)-
\int_{\partial U_{r}(s)}f(z)d_{z}^{c}\log\|d\vartheta_{z}\|^{2}_{\mathrm{res,can}}(z)\bigg)\xrightarrow[r\rightarrow 0]{}\\&\sum_{\substack{s\in \mathrm{Sing}(f)\\
s\not\in \mathcal{P}}}\frac{c_{f,s}}{2}\log\|d\vartheta_{z}\|^{2}_{\mathrm{res,can}}(s),\\
&\,\,\,\,\,\,\sum_{p\in\mathcal{P}}\bigg(\int_{\partial U_{r}(p)}\log\|d\vartheta_{z}\|^{2}_{\mathrm{res,can}}(z)d_{z}^{c}f(z)-
\int_{\partial U_{r}(p)}f(z)d_{z}^{c}\log\|d\vartheta_{z}\|^{2}_{\mathrm{res,can}}(z)\bigg) \xrightarrow[r\rightarrow 0]{}0,
\end{align*}
which completes the proof of the proposition. 
\end{proof}
\end{prop}
\begin{prop}\label{2.3prop2}
Let $f\in C_{\ell,\ell \ell}(\overline{X})$, then we have the equality of integrals 
\begin{align*}
-&\int_{U_{r_{0}}}\log \|d\vartheta_{z}\|^{2}_{\mathrm{res,hyp}}(z)d_{z}d_{z}^{c}f(z)= \frac{1}{2\pi}
\int_{U_{r_{0}}}f(z)\hathyp(z) + \\&\int_{U_{r_{0}}}f(z)\left(\int_{0}^{\infty}\del
\widehat{K}_{\mathrm{hyp}}(t;z)dt \right)\hathyp(z) +\sum_{\substack{s\in \mathrm{Sing}(f)\\s\not\in \mathcal{P}}}
\frac{c_{f,s}}{2}\log \|d\vartheta_{z}\|^{2}_{\mathrm{res,hyp}}(s)-\\
&\int_{\partial U_{r_{0}}}\log\|d\vartheta_{z}\|^{2}_{\mathrm{res,hyp}}(z)d_{z}^{c}f(z)+
\int_{\partial U_{r_{0}}}f(z)d_{z}^{c}\log\|d\vartheta_{z}\|^{2}_{\mathrm{res,hyp}}(z).
\end{align*}
\begin{proof}
From Corollary \ref{cor9}, we know that $\log \|d\vartheta_{z}\|^{2}_{\mathrm{res,hyp}}(z)$ is smooth on $X$ and is $\log\log$-singular at the cusps. So the proof of 
the proposition follows from equation \eqref{chernformhyp} and employing similar arguments as in Proposition \ref{2.3prop1}. 
\end{proof}
\end{prop}
\begin{prop}\label{2.3prop3}
Let $f\in C_{\ell,\ell \ell}(\overline{X})$, then we have the equality of integrals 
\begin{align}
-&\int_{U_{r_{0}}}\log\|d\vartheta_{z}\|^{2}_{\mathrm{res,hyp}}(z)
d_{z}d_{z}^{c}f(z)=\notag \\&
 2g\int_{U_{r_{0}}}f(z)\hatcan(z)-2\int_{U_{r_{0}}}f(z)\hatshyp(z) +\sum_{\substack{s\in
\mathrm{Sing}(f)\\s\not\in \mathcal{P}}}\frac{c_{f,s}}{2}\log \|d\vartheta_{z}\|^{2}_{\mathrm{res,hyp}}(s)- 
\notag\\&\int_{\partial U_{r_{0}}}\log\|d\vartheta_{z}\|^{2}_{\mathrm{res,hyp}}(z)d_{z}^{c}f(z)+
\int_{\partial U_{r_{0}}}f(z)d_{z}^{c}\log\|d\vartheta_{z}\|^{2}_{\mathrm{res,hyp}}(z).\label{2.3prop3eqn}
\end{align}
\begin{proof}
Subtracting equation \eqref{2.3prop1eqn} from the desired equality in \eqref{2.3prop3eqn}, it follows that 
for any $r> 0$, it suffices to prove that
\begin{align}
-&\int_{U_{r_{0},r}}\log\|d\vartheta_{z}\|^{2}_{\mathrm{res,hyp}}(z)d_{z}d_{z}^{c}f(z)+
\int_{U_{r_{0},r}}\log\|d\vartheta_{z}\|^{2}_{\mathrm{res,can}}(z)d_{z}d_{z}^{c}f(z)-\notag\\
&2\int_{U_{r_{0},r}}f(z)\hatcan(z) +2\int_{U_{r_{0},r}}f(z)\hatshyp(z) \xrightarrow[r\rightarrow 0]{}
\notag\\&\sum_{\substack{s\in \mathrm{Sing}(f)\\s\not\in \mathcal{P}}}\frac{c_{f,s}}{2}\big(\log
\|d\vartheta_{z}\|^{2}_{\mathrm{res,hyp}}(s)-\log\|d\vartheta_{z}\|^{2}_{\mathrm{res,can}}(s)\big)-\notag\\&
\int_{\partial U_{r_{0}}}\log\|d\vartheta_{z}\|^{2}_{\mathrm{res,hyp}}(z)d_{z}^{c}f(z)+
\int_{\partial U_{r_{0}}}f(z)d_{z}^{c}\log\|d\vartheta_{z}\|^{2}_{\mathrm{res,hyp}}(z)+\notag\\&
\int_{\partial U_{r_{0}}}\log\|d\vartheta_{z}\|^{2}_{\mathrm{res,can}}(z)d_{z}^{c}f(z)-\int_{\partial 
U_{r_{0}}}f(z)d_{z}^{c}\log\|d\vartheta_{z}\|^{2}_{\mathrm{res,can}}(z).\label{2.3prop3eqn1}
\end{align}
From Proposition \ref{prop7}, for any $r>0$ and $z\in U_{r_{0},r}$, we have
\begin{align}
&\can(z)-\shyp(z)= -d_{z}d_{z}^{c}\phi(z),\label{2.3prop3eqn2}\\&\log\|d\vartheta_{z}\|_{\mathrm{res,hyp}}^{2}
(z)- \log\|d\vartheta_{z}\|_{\mathrm{res,can}}^{2}(z) = \lim_{w\rightarrow z}\big(\ghyp(z,w)-\gcan(z,w)\big)=
2\phi(z).\label{2.3prop3eqn3}
\end{align}
Therefore, using the above two equations and Stokes's theorem, the left-hand side of limit 
\eqref{2.3prop3eqn1} simplifies to give
\begin{align*}
&2\int_{\partial U_{r}}\phi(z)d_{z}^{c}f(z)-2\int_{\partial U_{r}}
f(z)d_{z}^{c}\phi(z)-\\&\int_{\partial U_{r_{0}}}\log\|d\vartheta_{z}\|^{2}_{\mathrm{res,hyp}}(z)d_{z}^{c}f(z)
+\int_{\partial U_{r_{0}}}f(z)d_{z}^{c}\log\|d\vartheta_{z}\|^{2}_{\mathrm{res,hyp}}(z)+\\&
\int_{\partial U_{r_{0}}}\log\|d\vartheta_{z}\|^{2}_{\mathrm{res,can}}(z)d_{z}^{c}f(z)-
\int_{\partial U_{r_{0}}}f(z)d_{z}^{c}\log\|d\vartheta_{z}\|^{2}_{\mathrm{res,can}}(z).
\end{align*}
From Corollary \ref{cor8}, we know that $\phi(z)$ is smooth on $X$ and is $\log\log$-singular at the cusps. So 
employing similar arguments as in Proposition \ref{2.3prop1}, and using equation \eqref{2.3prop3eqn3}, we compute
\begin{align*}
&2\int_{\partial U_{r}}\phi(z)d_{z}^{c}f(z)-2\int_{\partial U_{r}}f(z)d_{z}^{c}\phi(z)\xrightarrow[r\rightarrow 0]{}
\sum_{\substack{s\in \mathrm{Sing}(f)\\s\not\in \mathcal{P}}}c_{f,s}\phi(s)=
\\& \sum_{\substack{s\in \mathrm{Sing}(f)\\s\not\in \mathcal{P}}}\frac{c_{f,s}}{2}\big(\log\|dz\|^{2}_{\mathrm{res,hyp}}(s)-
\log\|dz\|^{2}_{\mathrm{res,can}}(s)\big),
\end{align*}
which completes the proof of the proposition.
\end{proof}
\end{prop}
\begin{thm}\label{2.3thm4}
Let $f\in C_{\ell,\ell \ell}(\overline{X})$, then we have the equality of integrals 
\begin{align}
&g\int_{\overline{X}}f(z)\hatcan(z) =\notag \\&\bigg(\frac{1}{4\pi}+\frac{1}{\vx(X)}
\bigg)\int_{\overline{X}}f(z)\hathyp(z) + \frac{1}{2}\int_{\overline{X}}f(z)\bigg(\int_{0}^{\infty}\del\widehat{K}_{\mathrm{hyp}}(t;z)dt\bigg)\hathyp(z).\label{2.3thm4eqn}
\end{align}
\begin{proof}
From the equality of differential forms described in equation \eqref{keyidentity}, 
for any $f \in C_{\ell,\ell \ell}(\overline{X})$, 
we have the desired equality of integrals (\ref{2.3thm4eqn}) on the compact subset $Y_{r_{0}}$. 

For any $f \in C_{\ell,\ell \ell}(\overline{X})$, combining Propositions \ref{2.3prop2} and \ref{2.3prop3} 
proves the desired equality of integrals (\ref{2.3thm4eqn}) on $U_{r_{0}}$, and completes the proof of the 
theorem. 
\end{proof}
\end{thm}
\begin{cor}\label{2.3cor5}
Let $f\in C_{\ell,\ell \ell}(\overline{X})$, then we have the equality of integrals 
\begin{align*}
& g\int_{X}f(z)\can(z) = \\&\bigg(\frac{1}{4\pi}+\frac{1}{\vx(X)} \bigg)
\int_{X}f(z)\hyp(z) + \frac{1}{2}\int_{X}f(z)\bigg(\int_{0}^{\infty}
\del \khyp(t;z)dt \bigg)\hyp(z). 
\end{align*}
\begin{proof}
The proof follows from Theorem \ref{2.3thm4} and the fact that there are only finitely many cusps of $X$.  
\end{proof}
\end{cor}
\begin{rem}
Observe that for a fixed $w\in X$, as a function in the variable $z$, the hyperbolic Green's function 
$\ghyp(z,w)\in C_{\ell,\ell\ell}(\overline{X})$. Hence, combining Corollary \ref{2.3cor5} and Proposition 
\ref{prop7}, we find 
\begin{align}\label{phi(z)final}
&\phi(z)= \frac{1}{2g}\int_{X}\ghyp(z,\zeta)
\left(\int_{0}^{\infty}\del \khyp(t;\zeta)dt\right)\hyp(\zeta)-\notag\\&\frac{1}{8g^{2}}\int_{X}
\int_{X}\ghyp(\zeta,\xi)\bigg(\int_{0}^{\infty}\del \khyp(t;\zeta)dt\bigg)
\bigg(\int_{0}^{\infty}\del \khyp(t;\xi)dt\bigg)\hyp(\xi)\hyp(\zeta).
\end{align}
The above equation allows us to express the canonical Green's function $\gcan(z; w)$ in terms of 
expressions involving only the hyperbolic heat kernel $\khyp(t;z; w)$ and the hyperobolic metric $\hyp(z)$. 
Hence, in the upcoming article \cite{anil}, equation \eqref{phi(z)final} serves as a starting 
point for the extension of the bounds for the canonical Green's function $\gcan(z; w)$ from $\cite{jk}$ to 
noncompact hyperbolic Riemann orbisurfaces of finite volume. 

\vspace{0.2cm}
Furthermore, as stated before the key identity has been the most crucial tool in the work of J.~Jorgenson and 
J.~Kramer. We hope that the extended version of the key identity leads to the extension of their work.
\end{rem}
{\small{}}
\vspace{0.3cm}
{\small{
Department of Mathematics, IIT Hyderabad\\
Yeddumailaram, Hyderabad 502205, India\\ arya@iith.ac.in}}
\end{document}